\documentclass[11pt,twoside]{article}

\usepackage{a4wide}
\usepackage{amsfonts}
\usepackage{amssymb}
\usepackage{amsmath}
\usepackage{graphicx}

\newcommand{\ignore}[1]{}

\def\@begintheorem#1#2{\par\bgroup{\sc #1\ #2. }\it\ignorespaces}
\def\@opargbegintheorem#1#2#3{\par\bgroup{\sc #1\ #2\ (#3). } \it\ignorespaces}
\def\@endtheorem{\egroup}
\newtheorem{theorem}{Theorem}[section]
\newtheorem{corollary}[theorem]{Corollary}
\newtheorem{lemma}[theorem]{Lemma}
\newtheorem{example}[theorem]{Example}
\newtheorem{proposition}[theorem]{Proposition}
\newtheorem{definition}[theorem]{Definition}
\newcommand{\bt}[1]{\begin{theorem}\label{#1}}
\newcommand{\bc}[1]{\begin{corollary}\label{#1}}
\newcommand{\bl}[1]{\begin{lemma}\label{#1}}
\newcommand{\be}[1]{\begin{example}\label{#1}}
\newcommand{\bp}[1]{\begin{proposition}\label{#1}}
\newcommand{\ba}[1]{\begin{algorithm}\rm\label{#1}}
\newcommand{\bd}[1]{\begin{definition}\rm\label{#1}}{\normalsize }
\newcommand{\bpr}{\noindent {\em Proof. }}
\newcommand{\et}{\end{theorem}}
\newcommand{\ec}{\end{corollary}}
\newcommand{\el}{\end{lemma}}
\newcommand{\ee}{\end{example}}
\newcommand{\ep}{\end{proposition}}
\newcommand{\ed}{\end{definition}}
\newcommand{\epr}{{\ \vbox{\hrule\hbox{%
\vrule height1.3ex\hskip0.8ex\vrule}\hrule}}\\\par}

\def\R{\mathbb{R}}
\def\Z{\mathbb{Z}}
\def \lexmin {{\rm lexmin}}

\def\ra{\overline}
\def\la{\overline}
\def\sh{\overline}

\begin{document}

\title{\bf Shifted Matroid Optimization}

\author{
Asaf Levin
\thanks{\small Technion - Israel Institute of Technology, Haifa, Israel.
Email: levinas@ie.technion.ac.il}
\and
Shmuel Onn
\thanks{\small Technion - Israel Institute of Technology, Haifa, Israel.
Email: onn@ie.technion.ac.il}
}

\date{}

\maketitle

\begin{abstract}
We show that finding lexicographically minimal $n$ bases in a matroid can
be done in polynomial time in the oracle model. This follows from a more general result
that the shifted problem over a matroid can be solved in polynomial time as well.


\end{abstract}

\section{Introduction}
\label{introduction}

Let $G$ be a connected graph and let $n$ be a positive integer. Given $n$ spanning 
trees in $G$, an edge is {\em vulnerable} if it is used by all trees. 
We wish to find $n$ spanning trees with minimum number of vulnerable edges.
One motivation for this problem is as follows. We need to make a sensitive
broadcast over $G$. In the planning stage, $n$ trees are chosen and prepared.
Then, just prior to the actual broadcast, one of these trees is randomly chosen and used.
An adversary, trying to harm the broadcast and aware of the prepared trees
but not of the tree finally chosen, will try to harm a vulnerable edge,
used by all trees. So we protect each vulnerable edge with high cost,
and our goal is to choose $n$ spanning trees with minimum number of vulnerable edges.

Here we consider the following harder problem. For $k=1,\dots,n$,
call an edge {\em $k$-vulnerable} if it is used by at least $k$ of the $n$ trees.
We want to find $n$ {\em lexicographically minimal} trees,
that is, which first of all minimize the number of $n$-vulnerable edges,
then of $(n-1)$-vulnerable edges, and so on. More precisely, given $n$ trees,
define their {\em vulnerability vector} to be $f=(f_1,\dots,f_n)$ with $f_k$ the
number of $k$-vulnerable edges. Then $n$ trees with vulnerability vector $f$ are better
than $n$ trees with vulnerability vector $g$ if the last nonzero entry of $g-f$ is positive.
(We remark that this order is often used in the symbolic computation literature, where it is called
{\em reverse lexicographic}, but for brevity we will simply call it here {\em lexicographic}.)
As a byproduct of our results we show how to find in polynomial time $n$ lexicographically
minimal spanning trees, which in particular minimize the number of vulnerable edges.

This problem can be defined for any combinatorial optimization set as follows.
For matrix $x\in\R^{d\times n}$ let $x^j$ be its $j$-th column.
Define the {\em $n$-product of} a set $S\subseteq\R^d$ by
$$S^n\ :=\ \times_n S\ :=\ \{x\in\R^{d\times n}\ :\ x^j\in S\,,\ j=1,\dots,n \}\ .$$
Call two matrices $x,y\in\R^{d\times n}$ {\em equivalent} and write $x\sim y$ if
each row of $x$ is a permutation of the corresponding row of $y$.
The {\em shift} of matrix $x\in\R^{d\times n}$ is the unique matrix ${\la x}\in\R^{d\times n}$
which satisfies ${\la x}\sim x$ and $x^1\geq \cdots\geq x^n$, that is, the unique matrix
equivalent to $x$ with each row nonincreasing. Let $|x^j|:=\sum_{i=1}^d |x_{i,j}|$
and $|x|:=\sum_{j=1}^n|x^j|$ be the sums of absolute values of the components of $x^j$ and $x$,
respectively. The {\em vulnerability vector} of $x\in\{0,1\}^{d\times n}$ is
$$(|{\la x}^1|,\dots,|{\la x}^n|)\ .$$
We then have the following nonlinear combinatorial optimization problem.

\vskip.2cm\noindent{\bf Lexicographic Combinatorial Optimization.}
Given $S\subseteq\{0,1\}^d$ and $n$, solve
\begin{equation}\label{lexicographic_problem}
\lexmin\{(|{\la x}^1|,\dots,|{\la x}^n|)\ :\ x\in S^n\}\ .
\end{equation}

\vskip.2cm\noindent
The complexity of this problem depends on the presentation of the set $S$.
In \cite{KOS} it was shown that it is polynomial time solvable
for $S=\{z\in\{0,1\}^d\,:\,Az=b\}$ for any totally unimodular $A$ and any integer $b$.
Here we solve the problem for matroids.

\bt{matroid_theorem}
The lexicographic combinatorial optimization problem \eqref{lexicographic_problem}\,over the
bases of any matroid given by an independence oracle and any $n$ is polynomial time solvable.
\et
Our spanning tree problem is the special case with $S$ the set of indicators of spanning
trees in a given connected graph with $d$ edges, and hence is polynomial time solvable.

We note that Theorem \ref{matroid_theorem} provides a solution of a nonlinear optimization
problem over matroids, adding to available solutions of other nonlinear optimization
problems over matroids and independence systems in the literature,
see e.g. \cite{BLMORWW,BLOW,LOW,Onn} and the references therein.

We proceed as follows. In Section 2 we discuss the
{\em shifted combinatorial optimization problem} and its relation to the lexicographic
combinatorial optimization problem. In section 3 we solve the shifted problem
over matroids in Theorem \ref{matroid_theorem_2} and conclude Theorem \ref{matroid_theorem}.
In Section 4 we discuss matroid intersections, partially solve the shifted problem over
the intersection of two strongly base orderable matroids (which include gammoids)
in Theorem \ref{matroid_intersection_theorem}, and leave open the complexity
of the problem for the intersection of two arbitrary matroids.
We conclude in Section 5 with some final remarks about the polynomial time
solvability of the shifted and lexicographic problems over totally unimodular systems
from \cite{KOS}, and show that these problems are NP-hard over matchings already for $n=2$
and cubic graphs.

\section{Shifted Combinatorial Optimization}
\label{shifted}

Lexicographic combinatorial optimization can be reduced to the following problem.

\vskip.2cm\noindent{\bf Shifted Combinatorial Optimization.}
Given $S\subseteq\{0,1\}^d$ and $c\in\Z^{d\times n}$, solve
\begin{equation}\label{shifted_problem}
\max\{{\ra c}\,{\la x}\ :\ x\in S^n\}\ .
\end{equation}

The following lemma was shown in \cite{KOS}.
We include the proof for completeness.

\bl{reduction_1}\cite{KOS}
The Lexicographic Combinatorial Optimization problem \eqref{lexicographic_problem} can be reduced
in polynomial time to the Shifted Combinatorial Optimization problem \eqref{shifted_problem}.
\el
\bpr
Define the following $c\in\Z^{d\times n}$, and note that it satisfies ${\ra c}=c$,
$$c_{i,j}\ :=\ -(d+1)^{j-1}\,,\quad i=1,\dots,d\,,\quad j=1,\dots,n\ .$$
Consider any two vectors $x,y\in S^{n}$, and suppose that the vulnerability vector
$(|{\la x}^1|,\dots,|{\la x}^n|)$ of $x$ is lexicographically smaller than the
vulnerability vector $(|{\la y}^1|,\dots,|{\la y}^n|)$ of $y$. Let $r$ be the largest
index such that $|{\la x}^r|\neq |{\la y}^r|$. Then $|{\la y}^r|\geq |{\la x}^r|+1$.
We then have
\begin{eqnarray*}
{\ra c}\,{\la x}-{\ra c}\,{\la y}\ =\ c\,{\la x}-c\,{\la y}
&   =  & \sum_{j=1}^n(d+1)^{j-1}\left(|{\la y}^j|-|{\la x}^j|\right) \\
& \geq & \sum_{j<r}(d+1)^{j-1}\left(|{\la y}^j|-|{\la x}^j|\right)+(d+1)^{r-1} \\
& \geq & (d+1)^{r-1} -\sum_{j<r}d(d+1)^{j-1}\ >\ 0\ .
\end{eqnarray*}
Thus, an optimal solution $x$ for problem \eqref{shifted_problem}
is also optimal for problem \eqref{lexicographic_problem}.
\epr

We proceed to reduce the Shifted Combinatorial Optimization problem \eqref{shifted_problem}
in turn to two yet simpler auxiliary problems.
For a set of matrices $U\subseteq\{0,1\}^{d\times n}$ let $[U]$ be the set of
matrices which are equivalent to some matrix in $U$,
$$[U]\ :=\ \{x\in\{0,1\}^{d\times n}\ :\ \exists\ y\in U\,,\ x\sim y \}\ .$$
Consider the following two further algorithmic problems over a given $S\subseteq\{0,1\}^d$:

\begin{equation}\label{auxiliary_problem_1}
\mbox{{\bf Shuffling.} Given $c\in\Z^{d\times n}$, solve $\max\{cx\,:\,x\in[S^n]\}$.}
\end{equation}
\begin{equation}\label{auxiliary_problem_2}
\mbox{{\bf Fiber.} Given $x\in[S^n]$, find $y\in S^n$ such that $x\sim y$.}
\end{equation}

\bl{reduction_2}
The Shifted Combinatorial Optimization problem \eqref{shifted_problem} can be reduced in polynomial
time to the Shuffling and Fiber problems \eqref{auxiliary_problem_1} and \eqref{auxiliary_problem_2}.
\el
\bpr
First solve the Shuffling problem \eqref{auxiliary_problem_1} with profit matrix $\ra c$ and
let $x\in[S^n]$ be an optimal solution. Next solve the Fiber problem \eqref{auxiliary_problem_2}
for $x$ and find $y\in S^n$ such that $x\sim y$.
We claim that $y$ is optimal for the Shifted Combinatorial Optimization problem
\eqref{shifted_problem}. To prove this, we consider any $z$ which is feasible in
\eqref{shifted_problem}, and prove that the following inequality holds,
$${\ra c}\,{\la y}\ =\ {\ra c}\,{\la x}\ \geq\ {\ra c}\,x\ \geq\ {\ra c}\,{\la z}\ .$$
Indeed, the first equality follows since $x\sim y$ and therefore we have ${\la y}={\la x}$.
The middle inequality follows since ${\ra c}$ is nonincreasing. The last inequality follows since
$z\in S^n$ implies that $\la z\in[S^n]$ and hence $\la z$ is feasible in \eqref{auxiliary_problem_1}.
So $y$ is indeed an optimal solution for problem \eqref{shifted_problem}.
\epr

\section{Matroids}
\label{matroids}

Define the {\em $n$-union} of a set $S\subseteq\{0,1\}^E$ where $E$ is any finite set to be the set
$$\vee_n S\ :=\ \{x\in\{0,1\}^E\ :\ \exists x_1,\dots,x_n\in S\,,\ x=\sum_{k=1}^n x_k\}\ .$$

Call $S\subseteq\{0,1\}^E$ a matroid if it is the set of indicators
of independent sets of a matroid over $E$. We will have matroids over
$E=[d]:=\{1,\dots,d\}$ and $E=[d]\times[n]$.
The following facts about $n$-unions of matroids are  well known, see e.g. \cite{Sch}.

\bp{partition_matroid}
For any matroid $S$ and any $n$ we have that $\vee_n S$ is also a matroid.
Given an independence oracle for $S$, it is possible in polynomial time to realize an
independence oracle for $\vee_n S$, and if $x\in\vee_n S$,
to find $x_1,\dots,x_n\in S$ with $x=\sum_{k=1}^nx_k$.
\ep

Define the {\em $n$-lift} of a set of vectors $S\subseteq\{0,1\}^d$
to be the following set of matrices,
$$\uparrow_n\!S\ :=\ \{x\in\{0,1\}^{d\times n}\ :\ \sum_{j=1}^n x^j\in S\}\ .$$

We need two lemmas.
\bl{partition_1}
For any set $S\subseteq\{0,1\}^d$ and any $n$ we have that
$[S^n]=\vee_n \uparrow_n\!S$ in $\{0,1\}^{d\times n}$.
\el
\bpr
Consider $x\in [S^n]$. Then $x\sim y$ for some $y\in S^n$,
and thus $y^j\in S$ for $j=1,\dots,n$. Since $x\sim y$,
each row of $x$ is a permutation of the corresponding row of $y$.
Assume that the $i$-th row of $x$ is given by the permutation $\pi_i$
of the corresponding row of $y$.  That is, $x_{i,j}=y_{i,\pi_i(j)}$.
For $k=1,\dots,n$, we define a matrix $z_k\in \{0,1\}^{d\times n}$
whose column sum satisfies $\sum_{j=1}^n z_k^j=y^k$, by $(z_k)_{i,j}:=x_{i,j}$ if $\pi_i(j)=k$,
and otherwise $(z_k)_{i,j}:=0$.  Since $y^k\in S$, we conclude that $z_k\in \uparrow_n\!S$.
Since the supports of the $z_k$ are pairwise disjoint, we have
that $\sum_{k=1}^n z_k \in \vee_n \uparrow_n\!S$.
However, $\sum_{k=1}^n z_k=x$ by definition, and thus $x\in \vee_n \uparrow_n\!S$.

In the other direction, assume that $x\in \vee_n \uparrow_n\!S$.  Then there are
$x_1,\dots,x_n\in \uparrow_n\!S$ such that $x=\sum_{k=1}^n x_k$.
That is, there are $x_1,\dots,x_n \in \{0,1\}^{d\times n}$ such that for each
$k=1,\dots,n$, we have $\sum_{j=1}^n x_k^j \in S$, and $x=\sum_{k=1}^n x_k$.
Let $y$ be a matrix whose $k$-th column is $y^k=\sum_{j=1}^n x_k^j$.
Then $y^k\in S$, and therefore $y\in S^n$.  The matrices $x$ and $y$ are
$0-1$ matrices and $\sum_{j=1}^n x^j=\sum_{j=1}^n y^j$, and therefore $x\sim y$.
Thus, $x\in [S^n]$, as required.
\epr

\bl{uparrow_matroid}
Let $S\subseteq\{0,1\}^d$ be a matroid given by an independence oracle. Then we have:
\begin{enumerate}
\item
$\uparrow_n\!S\subseteq\{0,1\}^{d\times n}$ is also a matroid,
for which an independence oracle can be realized.
\item The rank of the matroid $\uparrow_n\!S$ equals to the rank of the matroid $S$.
\item The rank of the matroid $[S^n]=\vee_n \uparrow_n S$ equals $n$
times the rank of the matroid $S$.
\end{enumerate}
\el
\bpr
We begin with the first claim.  Let $z$ be the $0$ matrix (where $z_{i,j}=0$ for all $i,j$).
Then $\sum_{j=1}^n z^j$ is the zero vector that belongs to $S$ since $S$ is a matroid, and
therefore $z\in\uparrow_n\!S$.  Let $x\geq y$ be a pair of $0-1$ matrices such that
$x\in \uparrow_n\!S$.  Let $\hat{x}=\sum_{j=1}^n x^j$, and $\hat{y}=\sum_{j=1}^n y^j$.
Then $\hat{x}\geq \hat{y}$.  Since $x\in \uparrow_n\!S$, $\hat{x}\in S$, and because $S$ is
a matroid $\hat{y}\in S$, and thus $y\in \uparrow_n\!S$.  Last, assume that $x,y\in\uparrow_n\!S$
and $|x|>|y|$.  Let $\hat{x}=\sum_{j=1}^n x^j$ and $\hat{y}=\sum_{j=1}^n y^j$. Since $|x|>|y|$
we conclude that $|\hat{x}|>|\hat{y}|$.  Since $x,y\in\uparrow_n\!S$, $\hat{x},\hat{y}\in S$
and since $S$ is a matroid there is $i$ such that $\hat{x}_i=1$, $\hat{y}_i=0$ and changing
$\hat{y}_i$ to $1$ results in an indicating vector of an independent set of $S$.
Consider this value of $i$, and let $j$ be such that $x_{i,j}=1$. Then, by changing $y_{i,j}$
to $1$, we get a larger matrix whose column sum is the vector resulted from $\hat{y}$ by
changing its $i$-th component to $1$, and thus the new matrix is in $\uparrow_n\!S$.
Therefore, $\uparrow_n\!S$ is a matroid. To present an independence oracle of $\uparrow_n\!S$,
assume that we would like to test the independence of a $0-1$ matrix $x$, then we compute
$\sum_{j=1}^n x^j$ and test if the resulting vector is independent in $S$.

Next we show the second claim.  Let $x$ be an indicating matrix of an independent set
of the matroid $\uparrow_n\!S$, then its column sum is in $S$ and its rank in $S$ equals
to the rank of $x$ in $\uparrow_n\!S$.  In the other direction, let $\hat{x}\in S$.
Define a matrix $x$ whose first column is $\hat{x}$ and its all other columns are zero columns,
then $\sum_{j=1}^n x^j=\hat{x}\in S$, and therefore $x\in \uparrow_n\!S$. The rank of $x$
in $\uparrow_n\!S$ equals the rank of $\hat{x}$ in $S$, and thus the claim follows.

Last, consider the remaining claim.  Given a matroid $T$ it is always the case that the rank
of the ground set according to $\vee_n T$ is at most $n$ times the rank of the ground set
according to $T$.  Thus, to show the claim using the second claim, it suffices to show that given
a base $z$ of $S$, there is an independent set of $\vee_n \uparrow_n S$ whose rank
(in $\vee_n \uparrow_n S$) is $n$ times larger than the rank of $z$ (in $S$).
With a slight abuse of notation, assume that $z$ is the indicator of the base $z$.
Let $x$ be a matrix such that for $j=1,\dots,n$, $x^j=z$. Then $x\in S^n$, and therefore it is
an indicator of an independent set of $\vee_n \uparrow_n S$ whose rank is $|x|=n|z|$ that
equals $n$ times the rank of the base $z$ (according to the rank function of the matroid $S$).
\epr

\bt{matroid_theorem_2}
The shifted combinatorial optimization problem \eqref{shifted_problem} over any matroid given
by an independence oracle or over its set of bases and any $n$ is polynomial time solvable.
\et
\bpr
We begin with the problem over (the independent sets of) a matroid $S$.
Consider problem \eqref{auxiliary_problem_1} over $S$. Since $[S^n]=\vee_n \uparrow_n\!S$
by Lemma \ref{partition_1}, it follows from Proposition \ref{partition_matroid} and
Lemma \ref{uparrow_matroid} that $[S^n]$ is a matroid for which we can realize an independence
oracle. So the greedy algorithm over $[S^n]$ solves problem \eqref{auxiliary_problem_1} over $S$.
Next consider problem \eqref{auxiliary_problem_2} over $S$. Let $x\in[S^n]$ be given.
By Proposition \ref{partition_matroid} and Lemma \ref{partition_1} and Lemma \ref{uparrow_matroid}
again, we can find $x_1,\dots,x_n\in\uparrow_n\!S$ with $x=\sum_{k=1}^nx_k$.
Let $y$ be the matrix with $y^k=\sum_{j=1}^nx_k^j$ for $k=1,\dots,n$. Then
$$\sum_{j=1}^nx^j\ =\ \sum_{j=1}^n\sum_{k=1}^nx_k^j\ =\
\sum_{k=1}^n\sum_{j=1}^nx_k^j\ =\ \sum_{k=1}^ny^k\ .$$
This implies that $x\sim y$. Now, since $x_k\in\uparrow_n\!S$ and
$y^k=\sum_{j=1}^nx_k^j$, we have that $y^k\in S$ for all $k$.
So we can find $y\in S^n$ with $x\sim y$ and therefore solve problem \eqref{auxiliary_problem_2}
over $S$ as well. It now follows from Lemma \ref{reduction_2} that we can
indeed solve problem \eqref{shifted_problem} over $S$.

To solve problem \eqref{shifted_problem} over the set of bases of $S$ proceed as follows.
Define a new profit matrix $w$ by $w_{i,j}:=c_{i,j}+2|c|+1$ for all $i,j$.
Solve problem \eqref{auxiliary_problem_1} over (the independent sets of) $S$ with
profit $w$ by the greedy algorithm over $[S^n]$. Let $r$ be the rank of $S$
so that by Lemma \ref{uparrow_matroid} the rank of $[S^n]$ is $nr$.
Consider any $x,y,z\in[S^n]$ with $x,y$ bases and $z$ not a basis. Then
$$wx\ =\ cx+(2|c|+1)|x|\ =\ cx+(2|c|+1)nr\ \geq\ -|c|+(2|c|+1)nr\ ,$$
$$wz\ =\ cz+(2|c|+1)|z|\ \leq\ |c|+(2|c|+1)(nr-1)\ .$$
So $wx-wz>0$ and hence the $w$-profit of any basis $x$ is larger than that
of any independent set $z$ which is not a basis, so the algorithm will output a basis.
Also, $wx-wy=cx-cy$ and hence the $w$-profit of basis $x$ is larger
than that of basis $y$ if and only if the $c$-profit of $x$ is larger than that
of $y$, so the algorithm will output a basis $x$ of $[S^n]$ with largest $c$-profit.
Now solve problem \eqref{auxiliary_problem_2} and find $y\in S^n$ with $x\sim y$ as before.
Since $S$ has rank $r$, each column $y^k\in S$ satisfies $|y^k|\leq r$.
Since $nr=|x|=|y|=\sum_{k=1}^n|y^k|$ we must have in fact $|y^k|=r$ so that $y^k$ is a basis
of $S$ for all $k$. Therefore $y$ is an optimal solution for problem
\eqref{shifted_problem} over the bases.
\epr

\vskip.2cm\noindent{\em Proof of Theorem \ref{matroid_theorem}.}
This follows at once from Lemma \ref{reduction_1} and Theorem \ref{matroid_theorem_2}.
\epr

\section{Matroid Intersections}
\label{matroid intersections}

We next consider the Shifted Combinatorial Optimization problem over matroid intersections.
We can only provide a partial solution for the class of {\em strongly base orderable} matroids,
introduced in \cite{Bru}, which strictly includes the class of {\em gammoids},
see \cite[Section 42.6c]{Sch}.

Let $S\subseteq\{0,1\}^d$ be a matroid and let ${\cal B}$ be the set of subsets of $[d]$
which are supports of bases of $S$. Then $S$ is {\em strongly base orderable} if for every pair
$B_1,B_2\in{\cal B}$ there is a bijection $\pi:B_1\rightarrow B_2$ such that
for all $I\subseteq B_1$ we have $\pi(I)\cup(B_1\setminus I)\in {\cal B}$
(where $\pi(I)=\{ \pi(i): i\in I\}$).

We need the following elegant result of \cite{DM}, see also \cite[Theorem 42.13]{Sch}.
\bp{intersection}
For any two strongly base orderable matroids $S_1$ and $S_2$ we have
$$(\vee_n S_1)\cap(\vee_n S_2)\ =\ \vee_n(S_1\cap S_2)\ .$$
\ep

We also need the following lemma.
\bl{strongly_based_lemma}
If $S$ is a strongly base orderable matroid, then $\uparrow_n\!S$ is also a strongly
base orderable matroid. Further, for any $S_1,S_2\subseteq\{0,1\}^d$ we have that
$\uparrow_n\!S_1\cap \uparrow_n\!S_2=\uparrow_n\!(S_1\cap S_2)$.
\el
\bpr
We begin with the first part. Let $S$ be a strongly base orderable matroid.
Then by Lemma \ref{uparrow_matroid}, $\uparrow_n\!S$ is also a matroid.
Let $x,y\in\{ 0,1\}^{d\times n}$ whose supports $X,Y$ are two bases of $\uparrow_n\!S$.
By Lemma \ref{uparrow_matroid}, the supports $\hat{X}$ and $\hat{Y}$ of $\hat{x}=\sum_{j=1}^n x^j$
and $\hat{y}=\sum_{j=1}^n y^j$ respectively are two bases of $S$. Since $S$ is a
strongly base orderable matroid there is a bijection $\pi:\hat{X}\rightarrow \hat{Y}$
such that $\pi(\hat{I})\cup (\hat{X}\setminus \hat{I})$ is a base of $S$ for every
$\hat{I} \subseteq \hat{X}$. We define a bijection $\pi':X\rightarrow Y$ as follows.
For every $(i,j)\in X$, we let $\pi'(i,j)=(\pi(i),j')$ where $j'$ is the unique column
containing $1$ in the $\pi(i)$ row of $y$ (there is such a column as $\pi(i)\in \hat{Y}$).
A support of a $0-1$ matrix is a base in $\uparrow_n\!S$ if and only if its column sum is a
base in $S$. Thus, the required properties of $\pi'$ follow from these properties of $\pi$.

Next we prove the second claim. Let $x\in \uparrow_n\!S_1\cap \uparrow_n\!S_2$ and denote
its column sum by $\hat{x}=\sum_{j=1}^n x^j$. Since $x\in \uparrow_n\!S_1$, we have
$\hat{x}\in S_1$, and similarly since $x\in \uparrow_n\!S_2$, we have $\hat{x}\in S_2$.
Thus, $\hat{x}\in S_1\cap S_2$, and thus $x \in \uparrow_n\!(S_1\cap S_2)$.
In the other direction, let $x \in \uparrow_n\!(S_1\cap S_2)$ and denote
$\hat{x}=\sum_{j=1}^n x^j$. Then, $\hat{x} \in S_1\cap S_2$. Thus, $\hat{x} \in S_1$,
and therefore $x\in  \uparrow_n\!S_1$, and similarly $\hat{x} \in S_2$, and therefore
$x\in\uparrow_n\!S_2$.  Therefore, $x\in \uparrow_n\!S_1\cap \uparrow_n\!S_2$.
\epr

\bt{matroid_intersection_theorem}
The optimal objective function value of the shifted combinatorial optimization problem
\eqref{shifted_problem} over the intersection of any two strongly base orderable
matroids given by independence oracles and any $n$ can be computed in polynomial time.
\et
\bpr
Let $S_1,S_2\subseteq\{0,1\}^d$ be strongly base orderable matroids and let $S:=S_1\cap S_2$.
By Lemma \ref{strongly_based_lemma}, $\uparrow_n S_1$ and $\uparrow_n S_2$ are strongly base
orderable matroids. Applying Proposition \ref{intersection} to $\uparrow_n\!S_1$,
$\uparrow_n\!S_2$, and using Lemmas \ref{partition_1} and \ref{strongly_based_lemma},
we get $$[S_1^n]\cap[S_2^n]\ =\ (\vee_n \uparrow_n\!S_1)\cap(\vee_n \uparrow_n\!S_2)
\ =\ \vee_n(\uparrow_n\!S_1\cap \uparrow_n\!S_2)\ =\ \vee_n\uparrow_n\!(S_1\cap S_2)=[S^n]\ .$$
By Proposition \ref{partition_matroid} applied to the matroids $\uparrow_n\!S_1$ and
$\uparrow_n\!S_2$, and using Lemma \ref{partition_1}, we can realize in polynomial time independence
oracles for $[S_1^n]$ and $[S_2^n]$. So we can maximize a given profit matrix $c$ over $[S^n]$
by maximizing $c$ over the intersection of the matroids $[S_1^n]$ and $[S_2^n]$.
This solves problem \eqref{auxiliary_problem_1} over $S$ and gives the optimal value
of the shifted optimization problem \eqref{shifted_problem} over $S$ via
the equality $\max\{{\ra c}\,{\la x}\,:\,x\in S^n\}=\max\{cx\,:\,x\in [S^n]\}$.
\epr
We raise the following natural questions. Is it possible to actually find an optimal solution
(and not only the optimal value) for the shifted problem over the intersection of
strongly base orderable matroids in polynomial time? What is the complexity of
shifted combinatorial optimization over the intersection of two arbitrary matroids?

The answer to the first question is yes for transversal matroids since then the Fiber
problem \eqref{auxiliary_problem_2} over the intersection can also be solved.
This implies in particular that the shifted problem
over matchings in bipartite graphs is polynomial time solvable.
However, if the Fiber problem cannot be solved efficiently, that is given $x\in [S^n]$,
it is hard to find $y\in S^n$ such that $x\sim y$, then there are matrices $c$ for which
it is hard to solve the Shifted Combinatorial Optimization problem as well.
To see this, let $c$ be such that for all $i,j$ we have $c_{i,j}=1$ if $x_{i,j}=1$
and $c_{i,j}=-1$ if $x_{i,j}=0$. Then, the optimal profit will be $|x|$, and every optimal
solution $y$ to the Shifted Combinatorial Optimization problem must satisfy $y\in S^n$
and $y\sim x$, and thus solves the Fiber problem as well. Thus, the Fiber problem is
not harder than the Shifted Combinatorial Optimization problem.

\section{Remarks}
\label{remarks}

It was shown in \cite{KOS} that problem \eqref{shifted_problem} and hence
also problem \eqref{lexicographic_problem} are polynomial time solvable
for $S=\{z\in\{0,1\}^d\,:\,Az=b\}$ with $A$ totally unimodular and $b$ integer.
We briefly give a new interpretation of this by showing how the new auxiliary
problems \eqref{auxiliary_problem_1} and \eqref{auxiliary_problem_2} over such $S$
can also be efficiently solved, so Lemma \ref{reduction_2} can be applied.

\bp{unimodular}
For $S=\{z\in\{0,1\}^d:Az=b\}$ with $A$ totally unimodular and $b$ integer,
\begin{eqnarray}\label{unimodular_equation}
[S^n]\ =\ \{x\in\{0,1\}^{d\times n}\ :\ A\sum_{j=1}^nx^j=nb\}\ .
\end{eqnarray}
Moreover, problems \eqref{auxiliary_problem_1} and \eqref{auxiliary_problem_2}
hence \eqref{shifted_problem} and \eqref{lexicographic_problem}
over $S$ are polynomial time solvable.
\ep
\bpr
If $x$ is in the left hand side of \eqref{unimodular_equation} then there is a
$y\in S^n$ with $x\sim y$. Then $A\sum_{j=1}^nx^j=A\sum_{j=1}^ny^j=nb$ so $x$
is also in the right hand side. If $x$ is in the right hand side, then a decomposition
theorem of \cite{BT} implies that there is a $y\in S^n$ with $\sum_{j=1}^ny^j=\sum_{j=1}^nx^j$.
Then $x\sim y$ so $x$ is also in the left hand side of \eqref{unimodular_equation}.
Moreover, by an algorithmic version of \cite{KOS}
of this decomposition theorem, such a $y$ can be found in polynomial time,
solving problem \eqref{auxiliary_problem_2} over $S$. By the equality in
\eqref{unimodular_equation}, to solve problem \eqref{auxiliary_problem_1} over $S$ with a given
profit $c$, we can maximize $c$ over the right hand side of \eqref{unimodular_equation},
and this can be done in polynomial time by linear programming since the defining matrix
of this set is $[A,...,A]$, which is totally unimodular since $A$ is.
Lemmas \ref{reduction_1} and \ref{reduction_2} now imply that problems \eqref{shifted_problem}
and \eqref{lexicographic_problem} over $S$ can also be solved in polynomial time.
\epr
This in particular implies that lexicographic combinatorial optimization over $s-t$ paths
in digraphs or perfect matchings in bipartite graphs is polynomial time doable.

\vskip.2cm
The shifted combinatorial optimization problem over $S\subseteq\{0,1\}^d$ in the special
case $n=1$ is just the standard combinatorial optimization problem $\max\{cx:x\in S\}$
over $S$. So a solution of the former implies a solution of the latter. Is the converse true?
The above results show that it is true for matroids, certain matroid intersections, and unimodular
systems including $s-t$ paths in digraphs and perfect matchings in bipartite graphs.
Unfortunately, as we next show, the answer in general is negative.

\bp{matching}
Let $S\subseteq\{0,1\}^d$ be the set of perfect matchings in a cubic graph.
Then the lexicographic problem \eqref{lexicographic_problem} and shifted
problem \eqref{shifted_problem} over $S$ are NP-hard already for fixed $n=2$.
\ep
\bpr
By Lemma \ref{reduction_1} it is enough to show that solving
$\lexmin\{(|{\la x}^1|,|{\la x}^2|)\ :\ x\in S^2\}$ is NP-hard.
We claim that deciding if there is an $x\in S^2$ with $|{\la x}^2|=0$ is NP-complete.
Indeed, there is such an $x$ if and only if there are two edge disjoint perfect matchings
in the given graph. Now this holds if and only if there are three
pairwise edge disjoint perfect matchings in the given graph, since it is cubic.
This is equivalent to the existence of a $3$ edge coloring, which is NP-complete
to decide on cubic graphs \cite{Hol}.
\epr
We point out that the set of perfect matchings in a graph can be written in the
form $S=\{z\in\{0,1\}^d\,:\,Az=b\}$ with $A$ the vertex-edge incidence matrix of
the graph and $b$ the all $1$ vector in the vertex space. Thus, Proposition \ref{matching}
shows, in contrast with Proposition \ref{unimodular}, that if $A$ is not totally unimodular,
then the shifted problem is NP hard over such sets $S$, even if we can efficiently do
linear optimization over $S$, as is the case for perfect matchings.

\vskip.2cm
The shifted optimization problem $\max\{{\sh c}\,{\sh x}:x\in S^n\}$
can be defined for any set $S\subset\Z^d$, and its complexity
will depend on the structure and presentation of $S$.
Consider sets of bounded cardinality. Note that even with $|S|=2$,
the number of feasible solutions is $|S^n|=2^n$, so exhaustive search
is not polynomial. However, we do have the following simple statement.
\bp{fixed} For any fixed positive integer $s$, the shifted optimization problem
over any set $S\subset\Z^d$ satisfying $|S|\leq s$, any $n$,
and any $c\in Z^{d\times n}$, can be solved in polynomial time.
\ep
\bpr
Let $S=\{z^1,\dots,z^m\}$ be the list of elements of $S$, with $m\leq s$ and $z^i\in\Z^d$ for all
$i$. Given $x\in S^n$, let $n_i:=|\{k:x^k=z^i\}|$ for $i=1,\dots,m$. Let $y\in S^n$ be the matrix
with first $n_1$ columns equal to $z^1$, next $n_2$ columns equal to $z^2$, and so on, with last
$n_m$ columns equal to $z^m$. Then $y\sim x$ (in fact, $y$ is obtained from $x$ by applying
the same permutation to each row), and therefore ${\sh c}\,{\sh x}={\sh c}\,{\sh y}$.
So the objective value of $x$ can be computed from $n_1,\dots,n_m$. Now, any integers
$0\leq n_i\leq n$ with $\sum_{i=1}^m n_i=n$ give a feasible $y\in S^n$ with these counts.
So, it is enough to go over all such $n_1,\dots,n_m$, for each compute the corresponding $y$
and the value ${\sh c}\,{\sh y}$, and pick the best. The number of such tuples is at most
$(n+1)^{m-1}=O(n^{s-1})$, since $n_m=n-\sum_{i=1}^{m-1} n_i$ is determined by the others,
which is polynomial in $n$ for fixed $s$.
\epr


\begin{thebibliography}{}

\bibitem{BT}
Baum, S., Trotter, L.E., Jr.:
Integer rounding and polyhedral decomposition for totally unimodular systems.
Lecture Notes in Economical and Mathematical Systems 157:15--23 (1978)

\bibitem{BLMORWW}
Berstein, Y., Lee J., Maruri-Aguilar, H., Onn, S.,
Riccomagno, E., Weismantel, R., Wynn, H.:
Nonlinear matroid optimization and experimental design.
SIAM Journal on Discrete Mathematics 22:901--919 (2008)

\bibitem{BLOW}
Berstein, Y., Lee, J., Onn, S., Weismantel, R.:
Parametric nonlinear discrete optimization over
well-described sets and matroid intersections.
Mathematical Programming 124:233--253 (2010)

\bibitem{Bru}
Brualdi, R.A.: Common transversals and strong exchange systems.
Journal of Combinatorial Theory 8:307--329 (1970)

\bibitem{DM}
Davies, J., McDiarmid, C.:
Disjoint common transversals and exchange structures.
The Journal of the London Mathematical Society 14:55--62 (1976)

\bibitem{Hol}
Holyer, I.:
The NP-Completeness of Edge-Coloring.
SIAM Journal on Computing 10:718--720 (1981)

\bibitem{KOS}
Kaibel, V., Onn, S., Sarrabezolles, P.:
The unimodular intersection problem.
Operations Research Letters 43:592--594 (2015)

\bibitem{LOW}
Lee, J., Onn, S., Weismantel, R.:
Approximate nonlinear optimization over weighted independence systems.
SIAM Journal on Discrete Mathematics 23:1667--1681 (2009)

\bibitem{Onn}
Onn, S.: Nonlinear Discrete Optimization.
Zurich Lectures in Advanced Mathematics,
European Mathematical Society (2010),
available online at: {\tt http://ie.technion.ac.il/$\sim$onn/Book/NDO.pdf}

\bibitem{Sch}
Schrijver, A.: Combinatorial Optimization, Springer (2003)

\end{thebibliography}
\end{document}